\documentclass[11pt]{article}
\usepackage{amsmath,amssymb,amsthm,color}
\usepackage[bookmarks,colorlinks]{hyperref}
\usepackage{natbib}

\newtheorem{theorem}{Theorem}[section]
\newtheorem{corollary}[theorem]{Corollary}
\theoremstyle{definition}

\theoremstyle{remark}
\newtheorem{remark}[theorem]{Remark}
\numberwithin{equation}{section}
\newcommand{\R}{\ensuremath{\mathbb{R}}}

\newcommand{\abs}[1]{\left\vert#1\right\vert}
\newcommand{\set}[1]{\left\{#1\right\}}
\newcommand{\paren}[1]{\left(#1\right)}
\newcommand{\brac}[1]{\left[#1\right]}

\begin{document}

\centerline{\Large{\bf A note on second order linear functional equations }}

\centerline{}

\centerline{\Large{\bf in random normed spaces }}

\centerline{}

\centerline{Mongkhon Tuntapthai}

\centerline{}

\centerline{Department of Mathematics, Faculty of Science}

\centerline{Khon Kaen University, Khon Kaen, Thailand}

\centerline{}

\centerline{mongkhon@kku.ac.th}

\begin{abstract}
In this paper, we apply the publication of Joung (2009) to derive a stability result for for the second order linear functional equation:
\[ f(x) = pf(x-1)-qf(x-2) \ \text{ for all } x\in\mathbb R,\]
where $f$ is a mapping from $\mathbb R$ into the induced random space of any Banach space. By relaxing the lower bound assumption, we also generalize the result of Jung (2009) on arbitrary random normed spaces with the minimum $t$-norm. However, we need the monotonicity of the distribution in the lower bound assumption. By the properties of normal distributions, our main result can be applied.
\end{abstract}




\section{Introduction}

A well-known problem in the theory of functional equations is the approximation of any approximately linear function by a linear function, called the Hyers-Ulam stability. 
It is based on a question of \cite{ulam1960collection} asking for the approximation of homomorphism between metric groups. 
The affirmatively answered of the specific Ulam\rq s question for linear mapping between Banach spaces was proved by \cite{hyers1941stability}.
Subsequently, many authors have been extended the Hyers result to some type of equations such as the Cauchy, quadratic, cubic and Jensen type functional equations, and we refer to \cite{aoki1950stability, rassias1978stability, rassias1998stability, jun2002generalized, lee2007stability} for more details.
Recently, \cite{soon2009functional,brzdkek2010note,jung2011hyers} proved the Hyers-Ulam stability and investigated a solution mapping $f\colon \R\to X$ of the second order linear functional equation:
\begin{eqnarray}
	f(x) = pf(x-1) - qf(x-2)
	\text{ for all } x\in\R,
\label{eq1.1}
\end{eqnarray}
where $X$ is a Banach space and $p$, $q$ are real numbers such that $q\ne 0$ and $p^2 - 4q \ne 0$. 
As a special case of $p=1$ and $q=-1$, the equation (\ref{eq1.1}) is called the Fibonacci functional equation. 
Further results for the linear functional equations of higher order were published in \cite{brzdkek2008hyers,kim2013functional,jung2014linear}. 


Throughout this paper, the space of distribution functions is given by $\Delta^+$, 
that is the space of all mappings $F \colon \R \cup \{+\infty,-\infty\} \to [0,1]$
such that $F$ is left continuous non decreasing and $F(0)=0$, $F(+\infty)= 1$. 
Denote $D^+$ by a subset of $\Delta^+$ for which $\ell^- F(+\infty) =1$, 
where $\ell^{-1} F(x)$ is the left limit of the function $F$ at the point $x$, that is, $\ell^{-} F(x) = \lim_{t\to x^-} F(t)$.

The space $\Delta^+$ is partially ordered by the usual point-wise ordering of functions, i.e., $F \le G$ if and only if $F(x) \le G(x)$ for all $x$ in $\R$. 
The maximal element for $\Delta^+$ in this order is the distribution function $\varepsilon_0$ given by
\[
	\varepsilon_0 (t) =
	\begin{cases}
	0 & \text{if } t\le 0
	\\
	1 & \text{if } t> 0.
	\end{cases}
\]

A mapping $T \colon [0, 1] \times [0, 1] \to [0, 1]$ is a \emph{triangular norm} (briefly, a $t$-norm) if $T$ satisfies the following conditions:
\begin{enumerate}
	\item[(TN1)] $T$ is commutative and associative;
	\item[(TN2)] $T(a, 1) = a$ for all $a \in [0, 1]$;
	\item[(TN3)] $T(a, b) \le T(c, d)$ whenever $a \le c$ and $b \le d$ for all $a, b, c, d \in [0, 1]$.
\end{enumerate}
A $t$-norm $T$ can be extended (by associativity) to an $n$-array operation taking for $(a_1,a_2,\dots,a_n) \in [0,1]^n$ , the value $T(a_1,a_2,\dots,a_n)$ defined by
\[
	T_{i=1}^0 a_i = 1,
	\quad
	T_{i=1}^n a_i = T(T_{i=1}^{n-1} a_i, a_n) 
	= T(a_1,a_2,\dots,a_n)
\]
Typical examples of $t$-norms are $T_P (a, b) = ab$ (the product of two numbers) and $T_M (a, b) = \min(a, b)$ (the minimum of two numbers).

A \emph{random normed space} (briefly, RN-space) is a triple $(X,\mu, T)$, where $X$ is a vector space, $T$ is a $t$-norm, and $\mu$ is a mapping from $X$ into $D^+$ ($\mu(x)$ denoted by $\mu_x$) such that, the following conditions hold:
\begin{enumerate}
	\item[(RN1)] $\mu_x(t) = \varepsilon_0 (t)$ for all $t > 0$ if and only if $x = 0$;
	\item[(RN2)] $\mu_{\beta x} (t) = \mu_x ( \frac{t}{|\beta|} )$ for all $x \in X$, $\beta \ne 0$;
	\item[(RN3)] $\mu_{x+y} (t + s) \ge T(\mu_x (t), \mu_y (s))$ for all $x, y \in X$ and $t, s \ge 0$.
\end{enumerate}

Let $(X,\mu, T)$ be an RN-space.
A sequence $\set{x_n}$ in $X$ is said to be \emph{convergent} to $x$ in $X$ if, 
	for every $\varepsilon > 0$ and $\lambda > 0$, there exists a positive integer $N$ 
	such that $\mu_{x_n -x} (\varepsilon) > 1 -\lambda$ whenever $n \ge N$.
A sequence $\set{x_n}$ in $X$ is called \emph{Cauchy sequence} if, 
	for every $\varepsilon > 0$ and $\lambda > 0$, there exists a positive integer $N$ 
	such that $\mu_{x_n -x_m} (\varepsilon) > 1 -\lambda$ whenever $n \ge m \ge N$.
An RN-space $(X,\mu, T)$ is said to be \emph{complete} if and only if every Cauchy sequence in $X$ is convergent to a point in $X$.


A number of publications have been generalized the Hyers-Ulam stability for various kinds of functional equations in an RN-space such as the Cauchy, quadratic and cubic functional equations, see for examples \cite{baktash2008stability, gordji2009approximation, kenary2009stability, mihect2010stability, schin2010stability,  chugh2012stability, cho2013stability}

In this article, we focus on the Hyers-Ulam stability for (\ref{eq1.1}) in the class of functions $f\colon \R \to X$, where $X$ is a random normed space. 


\section{Main Theorem}

The stability result for the functional equation (\ref{eq1.1}) in Banach spaces was given by \cite{soon2009functional}. 
Now, we generalize and improve the result of \cite{soon2009functional} on RN-spaces with the minimum $t$-norm.

\begin{theorem}
\label{thm2.1}
Let $p$ and $q$ be two real numbers such that the quadratic equation $x^2-px+q=0$ has distinct real solutions $\alpha$ and $\beta$ with $0<|\beta |<|\alpha |<1$. 
Assume that a mapping $\varphi \colon \R \to D^+$ (where $\varphi_x$ denotes $\varphi(x)$) satisfies the conditions
\[
	\lim_{x\to -\infty} \varphi_x(t) = 1 
	\text{ for all } t>0 
\]
and
\[
	\varphi_x \ge \varphi_{x+1} 
	\text{ for all } x\in \R .
\]
Let $(X, \mu, T_M)$ be a complete RN-space. 
If a mapping $f \colon \R \to X$ satisfies that 
\begin{eqnarray}
	\mu_{f(x)-pf(x-1)+qf(x-2)} 
	\ge \varphi_x
	\text{ for all } x\in\R ,
\label{eq2.1}
\end{eqnarray}
then there exists $F \colon \R \to X$ a solution of the functional equation (\ref{eq1.1}) such that
\[
	\mu_{f(x)-F(x)} (t)
	\ge 
	\varphi_x \paren{ \gamma t } 
	\text{ for all } x\in\R 
	\text{ and } t>0
\]
where 
\[
	\gamma = \tfrac{ \abs{\alpha -\beta} \, \paren{1-|\alpha | } \, \paren{1-|\beta | } }
	{ |\alpha | + |\beta | -2|\alpha | |\beta | } .
\]
\end{theorem}


\proof
For all $n\in \mathbb N$, define a sequence of mappings $G_n \colon \R \to X$ by
\[
	G_n(x) = \alpha^n \set{f(x -n) - \beta f(x - n-1)} 
	\text{ for all } x \in \R.
\]
Since $\alpha + \beta = p$ and $\alpha \beta = q$, for all $x \in \R$ and $t>0$,
\begin{eqnarray*}
	&&\hskip-.75cm
	\mu_{ G_n(x) - G_{n+1}(x) } \paren{t|\alpha |^n}
	\notag\\
	&= &
	\mu_{\alpha^n \set{f(x -n) - \beta f(x - n-1) }
	-\alpha^{n+1} \set{f(x -n-1) - \beta f(x - n-2)} }
	\paren{t|\alpha|^n}
	\notag\\
	&= &
	\mu_{ f(x -n) - \beta f(x - n-1) 
	-\alpha \set{f(x -n-1) - \beta f(x - n-2)} } 
	\paren{t}
	\notag\\
	&= &
	\mu_{ f(x -n) - pf(x - n-1) +q f(x - n-2) } 
	\paren{t} .
\end{eqnarray*}
Substituting $x-n$ for $x$ to the assumption (\ref{eq2.1}), we get
\begin{equation}
	\mu_{ G_n(x) - G_{n+1}(x) } \paren{t|\alpha |^n}
	\; = \; \mu_{ f(x -n) - pf(x - n-1) +q f(x - n-2) } \paren{t} 
	\; \ge \; \varphi_{x-n} (t) .
\label{eq2.2}
\end{equation}
The last inequality implies that for all $n,m \in \mathbb N$,
\begin{eqnarray}
	&&\hskip-.75cm
	\mu_{G_{n}(x) - G_{n+m}(x)} 
	\textstyle 
	\paren{t \sum\limits_{k=n}^{n+m-1} |\alpha|^k}
	\notag\\
	&= &
	\mu_{ \set{G_{n}(x) - G_{n+1}(x)} + \set{G_{n+1}(x) - G_{n+2}(x)} + \dots \set{ G_{n+m-1}(x) - G_{n+m}(x)} } 
	\textstyle
	\paren{t \sum\limits_{k=n}^{n+m-1} |\alpha|^k}
	\notag\\
	&\ge &
	T_M \left(
	\mu_{G_n(x) -G_{n+1}(x)} \paren{t|\alpha|^n} ,
	\mu_{G_{n+1}(x) -G_{n+2}(x)} \paren{t|\alpha|^{n+1}} ,
	\right.
	\notag\\
	&&
	\phantom{ T_M \big\{ \, \, }
	\left. \dots , \mu_{G_{n+m-1}(x) -G_{n+m}(x)} \paren{t|\alpha|^{n+m-1}}  \right)
	\notag\\
	&\ge &
	T_M \paren{\varphi_{x-n} (t) ,\varphi_{x-n-1} (t) ,\dots,\varphi_{x-n-m-1} (t)}
	\notag\\
	&= &
	\varphi_{x-n} (t) .
\label{eq2.3}
\end{eqnarray}
By (\ref{eq2.3}) the condition that $\lim_{x\to -\infty} \varphi_{x}(t) = 1$, $\set{G_n(x)}_{n=1}^\infty$ is a Cauchy sequence in $X$. 
Since $X$ is complete, we can define a mapping $G \colon \R \to X$ by
\[
	G(x) = \lim_{n\to \infty} G_n(x)
	\text{ for all } x\in\R.
\]
In view of (\ref{eq2.2}), we set $n=0$ and take $m$ to the infinity to obtain 
\begin{eqnarray}
	\mu_{f(x) -\beta f(x-1) - G(x)} 
	\paren{ \tfrac{t}{1-|\alpha|} }
	\ge \varphi_x (t) 
	\text{ for all } t>0.
\label{eq2.4}
\end{eqnarray}
Also, we observe that for all $x\in\R$,
\begin{eqnarray}
	pG(x - 1) - qG(x - 2)
	&= &
	p \alpha^{-1} \lim_{n\to\infty} \alpha^{n+1} 
	\brac{ f(x - (n + 1)) - \beta f(x - (n + 1) - 1) }
	\notag\\
	&&
	- q \alpha^{-2} \lim_{n\to\infty} \alpha^{n+2} 
	\brac{ f(x - (n + 2)) - \beta f(x - (n + 2) - 1) }
	\notag\\
	&= &
	p \alpha^{-1} G(x) - q \alpha^{-2} G(x)
	\notag\\
	&= &
	G(x) .
\label{eq2.5}
\end{eqnarray}
%


Form above argument, it is easy to show that a mapping $H \colon \R \to X$ given by
\[
	H(x) = \lim_{n\to \infty} H_n(x)
	\text{ for all } x\in\R
\]
is well-defiend, where the mappings $H_n \colon \R \to X$ are defiend by
\[
	H_n(x) = \beta^n \set{f(x -n) - \alpha f(x - n-1)} 
	\text{ for all } n\in \mathbb N.
\]
Moreover, we can check that the mapping $H$ satisfies that 
\begin{eqnarray}
	\mu_{f(x) -\alpha f(x-1) - H(x)} \paren{ \tfrac{t}{1-|\beta|} }
	\ge \varphi_x (t) 
	\text{ for all } t>0
\label{eq2.6}
\end{eqnarray}
and it solves the linear functional equation (\ref{eq1.1}), i.e. 
\begin{eqnarray}
	pH(x - 1) - qH(x - 2)
	&= &	H(x) 	
	\text{ for all } x\in\R .
\label{eq2.7}
\end{eqnarray}



Finally, we assert that a mapping $F\colon\R\to X$ given by 
\[
	F(x) = \frac{\alpha}{\alpha-\beta} G(x) - \frac{\beta}{\alpha-\beta} H(x) 
	\text{ for all } x\in\R ,
\]
satisfies the requirements of our main theorem.
By (\ref{eq2.5}) and (\ref{eq2.7}), we have
\begin{eqnarray*}
	&&\hskip-.75cm
	pF(x - 1) - qF(x - 2)
	\notag\\
	&= &
	\frac{p\alpha}{\alpha-\beta} G(x-1) - \frac{p\beta}{\alpha-\beta} H(x-1)
	- \frac{q\alpha}{\alpha-\beta} G(x-2) + \frac{q\beta}{\alpha-\beta} H(x-2)
	\notag\\
	&= &
	\frac{\alpha}{\alpha-\beta} \set{ p G(x-1) - qG(x-2) } 
	- \frac{\beta}{\alpha-\beta} \set{ pH(x-1) -qH(x-2) }
	\notag\\
	&= &
	\frac{\alpha}{\alpha-\beta} G(x)
	- \frac{\beta}{\alpha-\beta} H(x)
	\notag\\
	&= &
	F(x) 
\end{eqnarray*}
for all $x\in\R$, whence the mapping $F$ solves the linear functional equation (\ref{eq1.1}).
Note that  for all $x\in\R$,
\begin{eqnarray*}
	&&\hskip-.75cm
	f(x) - F(x) 
	\notag\\
	&= &
	\frac{1}{\alpha-\beta} \brac{\paren{\alpha-\beta} f(x) - \set{\alpha G(x) - \beta H(x)} }
	\notag\\
	&= & 
	\frac{\alpha}{\alpha-\beta} \brac{ f(x) - \beta f(x-1) -G(x) } 
	- \frac{\beta}{\alpha-\beta} \brac{ f(x) - \alpha f(x-1) - H(x) } ,
\end{eqnarray*}
so that for all $t>0$,
\begin{eqnarray*}
	&&\hskip-.75cm
	\mu_{f(x) - F(x) } \paren{ \tfrac{ \paren{ |\alpha | + |\beta | -2|\alpha | |\beta | } t }{\abs{\alpha -\beta} \paren{1-|\alpha | } \paren{1-|\beta | } } } 
	\notag\\
	&= &
	\mu_{f(x) - F(x) } \paren{ \tfrac{ |\alpha |t}{\abs{\alpha -\beta} \paren{1-|\alpha | } } 
	+ \tfrac{ |\beta |t}{\abs{\alpha -\beta} \paren{1-|\beta | } } } 
	\notag\\
	&\ge &
	T_M \left( \mu_{ \frac{\alpha}{\alpha-\beta} \brac{ f(x) - \beta f(x-1) -G(x) } }
	\paren{ \tfrac{ |\alpha |t}{\abs{\alpha -\beta} \paren{1-|\alpha | } } } ,
	\right.
	\notag\\
	&&
	\phantom{T_M ( \, \, }
	\left. \mu_{\frac{-\beta}{\beta-\alpha} \brac{ f(x) - \alpha f(x-1) - H(x) } } 
	\paren{ \tfrac{ |\beta |t}{\abs{\alpha -\beta} \paren{1-|\beta | } } } 
	\right) 
	\notag\\
	&= &
	T_M \paren{ \mu_{ f(x) - \beta f(x-1) -G(x) }  \paren{ \tfrac{t}{1-|\alpha |} } ,
	\mu_{ f(x) - \alpha f(x-1) - H(x) }  \paren{ \tfrac{t}{1-|\beta |} } } 
	\notag\\
	&\ge &
	T_M \paren{ \varphi_x (t ), \varphi_x (t ) }
	\notag\\
	&= &
	\varphi_x (t ) 
\end{eqnarray*}
where we have used (\ref{eq2.4}) and (\ref{eq2.6}) in the first inequality.
Therefore,
\[
	\mu_{f(x) - F(x) } (t)
	\ge \varphi_x \paren{ \gamma t}  
	\text{ for all } x\in\R
	\text{ and } t>0 .
\]
\endproof


For each $x \in \R$, define $\varphi_x \colon \R\to D^+$ by
\[
	\varphi_x(t)
	:= \frac{1}{\sqrt{2\pi}} \int_{-\infty}^{t-x} e^{\tfrac{-u^2}{2}} \,du
	\text{ for all } t>0 .
\]
where $\varphi_x$ denotes $\varphi (x)$. 
By changing of variables, we have 
\[
	\varphi_x(t)
	= \frac{1}{\sqrt{2\pi}} \int_{-\infty}^t \exp \paren{-\tfrac{(u-x)^2}{2}} \,du 
	\text{ for all } t>0 .
\]
Recall that the mapping $\varphi_x$ is the cumulative distribution function (called a normal distribution) which has the location parameter equal to $x$ and the scale parameter equal to $1$.  
Because the properties of the location parameter $x$ coincides with the conditions for the mapping $\varphi_x$ in Theorem \ref{thm2.1}, we obtain the following statement:
\begin{corollary}
Let $(X, \mu, T_M)$ be a complete RN-space and 
$p$, $q$ real numbers such that the quadratic equation $x^2-px+q=0$ has distinct real solutions $\alpha$ and $\beta$ with $0<|\beta |<|\alpha |<1$. 
If a mapping $f \colon \R \to X$ satisfies that for all $x\in\R$ and $t>0$,
\[
	\mu_{f(x)-pf(x-1)+qf(x-2)} (t)
	\ge \frac{1}{\sqrt{2\pi}} \int_{-\infty}^{t-x} e^{\tfrac{-u^2}{2}} \,du ,
\]
then there exists $F \colon \R \to X$ a solution of the functional equation (\ref{eq1.1}) such that
\[
	\mu_{f(x)-F(x)} (t)
	\ge 	\frac{1}{\sqrt{2\pi}} \int_{-\infty}^{\gamma (t-x)} e^{\tfrac{-u^2}{2}} \,du
\]
for all $x\in\R$ and $t>0$, where
\[
	\gamma = \tfrac{ \abs{\alpha -\beta} \, \paren{1-|\alpha | } \, \paren{1-|\beta | } }
	{ |\alpha | + |\beta | -2|\alpha | |\beta | } .
\]
\end{corollary}

\begin{remark} 
In the proof of main theorem, the inequality (\ref{eq2.2}) forces the mapping $\varphi_x$ to have the location parameter equal $x$. From this condition, we can prove the result only in the case of $0<|\alpha |<1$ and $0<\beta |<1$, where
\[
	\alpha = \frac{p+\sqrt{p^2-4q}}{2}
	\; \text{ and } \;
	\beta = \frac{p-\sqrt{p^2-4q}}{2}
\]
for some real numbers $p$ and $q$ such that $q\ne 0$ and $p^2-4q > 0$.
Besides the special case of $\alpha$ and $\beta$, we cannot guarantee that the mapping in Theorem \ref{thm2.1} is unique.
\end{remark}






\end{document}